\pageno=1
\documentstyle{amsppt}
\magnification=1200
\topmatter
\title On nonatomic Banach lattices and Hardy spaces \endtitle
\author N.J. Kalton and P. Wojtaszczyk \endauthor
\address{N.J. Kalton, Department of Mathematics, University of
Missouri, Columbia Mo. 65211, USA} \endaddress
\address{P. Wojtaszczyk, Institut of Mathematics, Warsaw University, 00-913
Warszawa, ul. Banacha 2, Poland} \endaddress
\thanks{ The  the first author was supported by NSF-grant
DMS-8901636; the second author was supported by KBN-grant 2-1055-91-01.
This work was done while the second named author was visiting
University of Missouri-Columbia. He would like to express his gratitude to the
whole Abstract Analysis Group of this Department for their hospitality}
\endthanks


\abstract{ We are interested in the question when a Banach space $X$ with
an
unconditional basis is isomorphic (as a Banach space) to an
order-continuous nonatomic Banach lattice. We
show that this is the case if and only if $X$ is isomorphic as a Banach space
with $X(\ell_2)$. This and results of J. Bourgain are used to show that spaces
$H_1(\bold T^n)$ are not isomorphic to nonatomic Banach lattices. We also show
that tent spaces introduced in \cite{4} are isomorphic to Rad $H_1$.
}
\endabstract

\endtopmatter

\document
\heading 1. Introduction \endheading
There is a natural distinction  between  sequence spaces and
function spaces (or, between atomic and nonatomic Banach lattices) in
functional analysis.
 As an example, let us point
out the subtitles of two volumes of \cite{15} and \cite{16}. However,
many classical function
spaces (e.g. the spaces $L_p[0,1]$ for $1<p<\infty$ \cite{22} or \cite
{16}) have unconditional bases
and hence are isomorphic as Banach spaces to sequence spaces (atomic
Banach lattices).  On the other hand, $L_1[0,1]$ and has no unconditional
basis (\cite {22} or \cite{16}) and in the other direction the sequence
spaces
$\ell_p$ for $p\neq 2$ are not isomorphic to any nonatomic Banach lattice
\cite{1}.
In this note we discuss a general criterion for deciding whether a Banach
space with an unconditional basis (i.e. a sequence space) can be
isomorphic to a nonatomic Banach lattice (i.e. a function space).
Our main result (Theorem 2.4) gives a simple necessary and sufficent
condition
for an atomic Banach lattice $X$ to be isomorphic to an order-continuous
nonatomic
Banach lattice; of course if $X$ contains no copy of $c_0$ every Banach
lattice structure on $X$ is order-continuous.

Our main motivation is to study the Hardy space $H_1(\bold T).$
After the discovery that the space $H_1(\bold T)$ has an unconditional
basis
 \cite{17} it become natural to
investigate if $H_1(\bold T)$ is isomorphic to a nonatomic Banach  lattice.
Applying Theorem 2.4 to $H_1$ and using some previous results of Bourgain
\cite {2} and \cite{3} we show that $H_1$ is not isomorphic to any
nonatomic Banach lattice, and further more that $H_1(\bold T^n)$ is not
isomorphic to a nonatomic Banach lattice for any natural number $n.$

We conclude by showing that the space
 Rad $H_1$ or $H_1(\ell_2)$ is isomorphic to the tent spaces $T^1$
introduced by Coifman, Meyer and Stein \cite{4}.

\heading  2. Lattices with unconditional bases\endheading

Our terminology about Banach lattices will agree with \cite {16}; we also
refer the reader to \cite{9} and \cite{10} for the isomorphic theory of
nonatomic Banach lattices.

A (real) Banach lattice $X$ is called {\sl order
continuous} if every order-bounded increasing sequence of positive
elements is norm convergent.  Any Banach lattice not containing $c_0$ is
automatically order-continuous.

For any order-continuous Banach lattice $X$ we can define an associated
Banach lattice
$ X(\ell_2)$ (using the Krivine calculus \cite {16} pp. 40-42) as the
space
of sequences $(x_n)_{n=1}^{\infty}$ in $X$ such that  $(\sum_{k=1}^n
|x_k|^2)^{1/2}$ is order-bounded (and hence is a convergent sequence) in
$X$.  $X(\ell_2)$ becomes an order-continuous Banach lattice when normed
by $\|(x_n)\|=\|(\sum_{n=1}^{\infty}|x_n|^2)^{1/2}\|.$

If $X$ has nontrivial cotype then $X(\ell_2)$ is naturally isomorphic to
the space Rad $X$ which is the subspace of $L_2([0,1];X)$ of functions of
the form $\sum_{n=1}^{\infty}x_nr_n$ where $(r_n)$ is the sequence of
Rademacher functions.  The space Rad $X$ is clearly an isomorphic
invariant of $X$, and so if two Banach lattices $X$ and $Y$ with
nontrivial cotype are isomorphic it follows easily that $X(\ell_2)$ and
$Y(\ell_2)$ are isomorphic.  However, this result holds in general by a
result of Krivine \cite{13} or \cite {16} Theorem 1.f.14.

\proclaim{Theorem 2.1} If $X$, $Y$ are order-continuous Banach lattices
and
$T:X\longrightarrow Y$ is a bounded linear operator, then  if $(x_n)\in
X(\ell_2)$ we have $(Tx_n)\in Y(\ell_2)$ and
$$ \|(T(x_n))\|_{Y(\ell_2)} \leq K_G \|T\| \|(x_n)\|_{X(\ell_2)}.$$\endproclaim
Here, as usual, $K_G$ denotes the Grothendieck constant.
\demo{Proof}Essentially this is Krivine's theorem, but we do need to show
that if $(x_n)\in X(\ell_2)$ then $(Tx_n)\in Y(\ell_2).$  To see this we
show that $(\sum_{k=1}^n|Tx_k|^2)^{1/2}$ is norm-Cauchy.  In fact if
$m>n$ then
$$ \align \|(\sum_{k=1}^m|Tx_k|^2)^{1/2}-(\sum_{k=1}^n|Tx_k|^2)^{1/2}\|_Y
&\le \|(\sum_{k=n+1}^m|Tx_k|^2)^{1/2}\|_Y\\
&\le K_G\|T\|\|(\sum_{k=n+1}^m|x_k|^2)^{1/2}\|_X\\
&\le K_G\|T\|\|(\sum_{k=n+1}^{\infty}|x_k|^2)^{1/2}\|_X
\endalign
$$
which converges to zero as $n\to\infty$ by the order-continuity of $X.$
\qed \enddemo

\proclaim{Corollary 2.2} If two order-continuous Banach lattices $X$ and
$Y$ are isomorphic as
Banach spaces, then $X(\ell_2)$ and $Y(\ell_2)$ are isomorphic as Banach
spaces.\endproclaim

If
$X$ is a separable order-continuous nonatomic Banach lattice then
$X$ can be
represented as (i.e. is linearily and order isomorphic with) a
K\"othe function space  on $[0,1]$ in such a way
that
$L_\infty[0,1]\subset X
\subset L_1[0,1]$ and inclusions are continuous.  It will then follow
that
 $L_\infty$ is dense in $X$, and
the dual of $X$ can be represented as a space of functions, namely
$X^*= \{ f\in L_1 \ : \ \int |fg|\, dt <\infty \ \text{
for every } \ g\in X \}.$

Now we are ready to state our main result.  Let us observe that for
re-arrangement invariant function spaces on $[0,1]$ this result was
proved in \cite{9} (cf. also \cite{16} 2.d) by a quite
different technique.
\proclaim{Theorem 2.3} Let $X$ be an order continuous, nonatomic Banach
lattice with an unconditional basis. Then $X$ is isomorphic as a Banach space
to $X(\ell_2)$. \endproclaim

\demo{Proof} We will represent $X$ as a K\"othe function space on
$[0,1]$ as described above.  Suppose
$(\phi_n)_{n=1}^{\infty}$ is a normalized
 unconditional basis of $X.$  Then there is an order-continuous atomic
Banach lattice
$Y$ which we identify as a sequence space
 and operators $U:X\longrightarrow Y$ and $V:Y
\longrightarrow X$ such that $UV=I_Y$, $VU=I_X$ and $U(\phi_n)=e_n$ for
$n=1,2,\dots$, where $e_n$ denotes the canonical basis vectors in $Y.$
We can regard $Y^*$ as a space of sequences and further suppose that
 $\|e_n\|_{Y^*}=\|e_n\|_Y=1.$  We will identify $Y(\ell_2)$ as a space of
double sequences with canonical basis $(e_{mn})_{m,n=1}^{\infty};$ thus
for any finitely nonzero sequence we have $\|\sum
a_{mn}e_{mn}\|_{Y(\ell_2)}=
\|\sum_m(\sum_n|a_{mn}|^2)^{1/2}e_m\|_Y.$

Let $r_n$ denote the Rademacher functions and for each fixed $f\in X$
note that $(r_nf)$ converges weakly to zero, since for $g\in X^*$ we have
$\lim_{n\to\infty}\int r_nfg\,dt=0.$
In particular we have for each $m\in\bold N$ that $(r_n\phi_m)$ converges
weakly to zero.  It follows by a standard gliding hump technique that if
$\eta=(2\|U\|\|V\|)^{-1}$ then we can find for each $(m,n)\in \bold N^2$
an integer $k(m,n)$  and disjoint subsets $(A_{mn})$ of $\bold N$ so that
$\|U(\phi_mr_{k(m,n)})\chi_{A_{mn}} - U(\phi_mr_{k(m,n)})\|_Y\le\eta.$

Identifying $Y^*$ as  a sequence space, we let $\psi_m=U^*(e_m)$ and then
define $v_{m,n}= \chi_{A_{mn}}U(\phi_mr_{k(m,n)})\in Y$ and
$v^*_{m,n}=\chi_{A_{mn}}V^*(\psi_mr_{k(m,n)})\in Y^*.$
Now suppose $(a_{mn})$ is a finitely nonzero double sequence.
Then
$$
\align
\| \sum_{m,n}a_{mn}v_{mn}\|_Y &\le
\|(\sum_{m,n}|a_{mn}|^2|U(\phi_mr_{k(m,n)})|^2)^{1/2}\|_Y \\
&\le K_G\|U\| \|(\sum_{m,n}|a_{mn}|^2|\phi_mr_{k(m,n}|^2)^{1/2}\|_X\\
&= K_G\|U\| \|(\sum_m(\sum_n |a_{mn}|^2)|\phi_m|^2)^{1/2}\|_X\\
&= K_G\|U\| \|(\sum_m(\sum_n |a_{mn}|^2)|Ve_m|^2)^{1/2}\|_X\\
&\le K_G^2\|U\|\|V\|
\|(\sum_m(\sum_n|a_{mn}|^2)^{1/2}|e_m|^2)^{1/2}\|_Y\\
&= K_G^2\|U\|\|V\| \|\sum_{m,n}a_{mn}e_{mn}\|_{Y(\ell_2)}.
\endalign
$$
  Here we have used
Krivine's theorem twice.
It follows that we can define a linear operator $S:Y(\ell_2)\to Y$ by
$Se_{mn}=v_{mn}$ and then $\|S\|\le K_G^2\|U\|\|V\|.$

Similar calculations yield that for any finitely nonzero double sequence
$(b_{mn})$ we have:
$$ \|\sum_{m,n}b_{mn}v_{mn}^*\|_{Y^*} \le K_G^2\|U\|\|V\|
\|\sum_m(\sum_n|b_{mn}|^2)^{1/2}e_{m}\|_{Y^*}.$$
Suppose then $y\in Y$  and set $a_{mn}=\langle y,v_{mn}^*\rangle.$
Let $F$ be a finite subset of $\bold N^2.$  Let
$\alpha_m=(\sum_n\chi_F(m,n)|a_{mn}|^2)^{1/2}$ and suppose the finitely
nonzero sequence $(\beta_m)$ is chosen so that $\|\sum
\beta_me_m\|_{Y^*}=1$ and $\sum \beta_m\alpha_m=\|\sum \alpha_me_m\|_Y.$
Then, with the convention that $0/0=0,$
$$
\align
\|\sum_{(m,n)\in F}a_{mn}e_{mn}\|_{Y(\ell_2)} &= \sum_m\beta_m\alpha_m\\
 &= \sum_{(m,n)\in F}\beta_m\alpha_m^{-1}|a_{mn}|^2\\
 &= \langle y, \sum_{(m,n)\in F}\beta_m\alpha_m^{-1}
a_{mn}v_{mn}^*\rangle\\
&\le \|y\|_Y
\|\sum_{(m,n)\in F}\beta_m\alpha_m^{-1}a_{mn}v^*_{mn}\|_{Y^*}\\
&\le K_G^2 \|U\|\|V\| \|y\|_Y.
\endalign
$$
Thus for each $F$ the map $T_F:Y\to Y(\ell_2)$ given by
$T_Fy=\sum_{(m,n)\in F}\langle y,v_{mn}^*\rangle e_{mn}$ has norm at most
$K_G^2\|U\|\|V\|.$  More generally, we have $\|T_Fy\| \le
K_G^2\|U\|\|V\|\|\chi_{A_F}y\|$ where $A_F=\cup_{(m,n)\in F}A_{mn}.$

 It follows  that for each $y\in Y$ the series $\sum_{m,n}
\langle y,v_{mn}^*\rangle e_{mn}$ converges (unconditionally) in
$Y(\ell_2).$ We can thus define an operator
$T:Y
\to Y(\ell_2)$ by $Ty=\sum_{m,n}\langle y,v_{mn}^*\rangle e_{mn}$ and
$\|T\|\le K_G^2\|U\|\|V\|.$

Now notice that $TS(e_{mn})=c_{mn}e_{mn}$ where $c_{mn}=\langle
v_{mn},v_{mn}^*\rangle.$  But:
$$
\align
 \langle v_{mn},v_{mn}^*\rangle &=\langle v_{mn},
V^*\psi_mr_{k(m,n)}\rangle\\
 &\ge \langle
U(\phi_mr_{k(m,n)}),V^*(\psi_mr_{k(m,n)})\rangle-\eta\|V\|\|\psi_m\|_{X^*}\\
 &=\langle \phi_m,\psi_m\rangle -\eta\|V\|\|\psi_m\|_{X^*}\\
 &\ge 1- \eta\|V\|\|U\| \ge 1/2.
 \endalign
$$
Thus $TS$ is invertible and so it follows that $Y(\ell_2)$ is isomorphic
to a complemented subspace of $Y.$  It then follows from the Pelczynski
decomposition technique that $Y\sim Y(\ell_2);$  more precisely
$Y\sim Y(\ell_2)\oplus W$ for some $W$ and so $Y\sim Y(\ell_2)\oplus
(Y(\ell_2)\oplus W) \sim Y(\ell_2)\oplus Y\sim Y(\ell_2).$
\qed \enddemo
\demo{Remark}  The order continuity of the Banach lattice $X$ is
essential.
In \cite{14} a non-atomic Banach lattice $X$ (actually an M-space) was
constructed which is isomorphic to $c_0$. In particular $X$ has an
unconditional basis but is not isomorphic to $X(\ell_2)$.  \enddemo

\proclaim{Theorem 2.4}Let $Y$ be a Banach space with an unconditional
basis.  Then $Y$ is isomorphic to an order-continuous nonatomic Banach
lattice if and only if $Y\sim Y(\ell_2).$
\endproclaim
\demo{Remark}Here again we regard $Y$ as an order-continuous Banach
lattice.
\enddemo
 \demo{Proof}One direction follows immediately from Theorem 2.3 and
Corollary 2.2.  For the other direction, it is only necessary to show
that if
$Y\sim Y(\ell_2)$ then
$Y$ is isomorphic to order-continuous nonatomic
Banach lattice.
To this we introduce the space $Y(L_2)$; this is the space of
sequences of functions $(f_n)$ in $L_2[0,1]$ such that
$\sum \|f_n\|_2 e_n$ converges in $Y$.  We set
$\|(f_n)\|_{Y(L_2)}=\|\sum\|f_n\|_2e_n\|_Y.$  It is clear that $Y(L_2)$
is an order-continuous Banach lattice.  Now if $(g_n)$ is an orthonormal
basis of $L_2$ we define $W:Y(\ell_2)\to Y(L_2)$ by
$W(\sum_{m,n}a_{mn}e_{mn}) = (\sum_n a_{mn}g_n)_{m=1}^{\infty}$
and it is easy to see that $W$ is an isometric isomorphism.\qed\enddemo

\proclaim{Proposition 2.5} If $X$ is a non-atomic order continuous
Banach
lattice with unconditional basis, then $X\sim X\oplus X$ and $X\sim X\oplus R$.
\endproclaim
\demo{Proof} Both facts follow from Theorem 2.3. \qed \enddemo

Note that for spaces with unconditional basis both properties do not hold in
general (see \cite{5} and \cite{6})
\proclaim{Proposition 2.6} Let $X$ be an order continuous
non-atomic Banach
lattice with an unconditional basis and let $Y$ be a complemented subspace of
$X$. Assume that $Y$ contains a complemented subspace isomorphic to $X$. Then
$X\sim Y$.\endproclaim
\demo{Proof} The proof  is a
repetition of the proof of Proposition 2.d.5. of \cite{16}. \qed \enddemo

\heading  3. Hardy spaces \endheading

We recall that $H_1(\bold T^n)$ is defined to be the space of boundary
values of functions $f$ holomorphic in the unit disk $\bold D$
and such that
$$\sup_{0<r<1} \int_{\bold T^n} |f(re^{it_1},re^{it_2},\dots,re^{it_n})|\, dt_1
dt_2\dots dt_n <\infty.$$ The basic theory of such spaces is explained in
\cite{18}.

Let us consider first the case $n=1$.  Then $\Re H_1$ is defined be the
space of real functions $f\in L_1(\bold T)$ such that for some $F\in
H_1(\bold T)$ we have $\Re F=f.$  $\Re H_1$ is normed by $\|f\|_1+
 \min\{ \|F\|_{H_1}: \Re F=f\}.$  Then $H_1$ is isomorphic to
the complexification of $\Re H_1$, and further when considered as a real
space is isomorphic to $\Re H_1.$  Further it was shown in \cite{17}
that
 $\Re H_1$ has an unconditional basis and is isomorphic a space of
martingales $H_1(\delta).$  To define the space $H_1(\delta)$ let
$(h_n)_{n\ge 1}$
be the usual enumeration of the Haar functions on $I=[0,1]$
 normalized so that $\|h_n\|_{\infty}=1.$  Then suppose $f\in
L_1$ is of the form
$f=\sum a_nh_n$.  We define $\|f\|_{H_1(\delta)}=\int(\sum_n
|a_n|^2h_n^2)^{1/2}dt$ and
$H_1(\delta)=\{f:\|f\|_{H_1(\delta)}<\infty\}.$

These considerations can be extended to the case $n>1.$  In a similar
way, $H_1(\bold T^n)$ is isomorphic to the complexification of, and is
also real-isomorphic to, a martingale space $H_1(\delta^n).$  Here we
define for $\alpha\in \Cal M=\bold N^n$ the function $h_{\alpha}\in
L_1(I^n)$ by $h_{\alpha}(t_1,\ldots,t_n)=\prod h_{\alpha_k}(t_k).$
Then $H_1(\delta^n)$ consists of all $f=\sum_{\alpha\in\Cal
M}a_{\alpha}h_{\alpha}$ such that $\|f\|_{H_1(\delta^n)}=\int (\sum
|a_{\alpha}|^2h_{\alpha}^2)^{1/2}dt<\infty.$

It is clear from the definition that the system $(h_{\alpha})_{\alpha\in
\Cal M}$ is an unconditional basis of $H_1(\delta^n).$  We can thus
define a space $H_1(\delta^n,\ell_2)=H(\delta^n)(\ell_2)$ as in Section
1; since
$H_1(\delta^n)$ has cotype two, this space is isomorphic to Rad
$H_1(\delta^n).$  The following theorem is due to Bourgain \cite{2}:

\proclaim{Theorem 3.1} $H_1(\delta,\ell_2)$ is not isomorphic to a
complemented subspace of $H_1(\delta).$\endproclaim

In a subsequent paper \cite{3} Bourgain implicitly extended this result
to
higher dimensions.

\proclaim {Theorem 3.2} For every $n=1,2,\dots$ the space
$H_1(\delta^n,\ell_2)$ is not isomorphic to any complemented subspace of
$H_1(\delta^n)$. \endproclaim
\demo{Sketch of proof} For $n=1$ this Theorem is proved in detail in \cite{2}.
The subsequent paper  \cite{3} states only the weaker
fact that $H_1(\delta^n)$ is not isomorphic to $H_1(\delta^{n+1})$. His proof
however gives the above Theorem as well. All that is needed is to change in
Section 3 of \cite{3} condition (m+1) and Lemma 4. Before we formulate the
appropriate condition we need some further notation. By $BMO(\delta^n)$ we will
denote the dual of $H_1(\delta^n)$ and by $BMO(\delta^n, \ell_2)$ we will
denote the dual of $H_1(\delta^n, \ell_2)$. The space $H_1(\delta^n,\ell_2)$
has an unconditional basis given by $(h_\alpha \otimes e_k)_{\alpha \in \Cal
M,\, k\in \bold N}$. In our notation from Section 2 $h_\alpha \otimes e_k$
is
a sequence of $H_1(\delta^n)$-functions which consists of zero functions except
at the $k$-th place where there is $h_\alpha$. The same element can be treated
as an element of the dual space.
Note that the natural duality gives
$$\langle h_\alpha \otimes e_k, h_{\alpha^\prime} \otimes e_{k^\prime} \rangle
=\cases \int_{I^n} |h_\alpha |, & \text{when $\alpha =\alpha^\prime $ and
$k=k^\prime$} \\ 0, &\text{otherwise.} \endcases $$
 Now we are ready to state the new condition (m+1):
\item{} Let $\Phi : H_1(\delta^n, \ell_2) \longrightarrow H_1(\delta^n)$
and $\Phi^{\times}: BMO(\delta^n, \ell_2) \longrightarrow BMO(\delta^n)$
be bounded
linear operators (note that $\Phi^{\times}$ is {\sl not} the adjoint of
$\Phi$). Then
for every $\varepsilon>0$ there exists a set $A\subset \Cal M$ such that
$\sum_{\alpha \in A} |h_\alpha| =1$ and integers $k_\alpha$ for $\alpha \in A$
such that
$$\sum_{\alpha \in A} \int_{I^n} | \Phi(h_\alpha \otimes e_{k_\alpha})|\cdot |
\Phi^{\times} (h_\alpha \otimes e_{k_\alpha})| <\varepsilon. $$

With this condition one can repeat the proof from \cite{3} and obtain the
Theorem. \qed  \enddemo

\proclaim{Corollary 3.3} We have the following
\define\cs{\overset c\to \subset}
$$\ell_2 \cs H_1(\delta) \cs H_1(\delta,\ell_2) \cs H_1(\delta^2)\cs
H_1(\delta^2, \ell_2) \cs \dots$$
where $X \cs Y$ means that $X$ is isomorphic to a complemented subspace of $Y$
but $Y$ is {\bf not} isomorphic to a complemented subspace of $X$. \endproclaim
\demo{Proof} It is well known and easy to check that the map $h_\alpha \otimes
e_k \mapsto h_\alpha (t_1,\dots ,t_n) \cdot r_k(t_{n+1})$ where $r_k$ is the
$k$-th Rademacher function gives the desired complemented embedding. That no
smaller space is isomorphic to a complemented subspace of a bigger one is the
above theorem of Bourgain. \qed \enddemo

\proclaim{Corollary 3.4} The spaces $H_1(\delta^n)$ is not isomorphic
to a nonatomic Banach lattice for
$n=1,2,\dots.$
The spaces $H_1(\delta^n,
\ell_2)$ are each isomorphic to a nonatomic Banach lattice.
\endproclaim
\demo{Proof} The first claim follows directly from Theorem 3.1, 3.2 and
Theorem
2.3. We only have to observe that (since $H_1(\delta^n)$ does not contain any
subspace isomorphic to $c_0$ and indeed has cotype two) any Banach
lattice isomorphic as a Banach space
to $H_1(\delta^n)$ is order continuous (see Theorem 1.c.4 of \cite {16}).  The
second claim follows from Corollary 2.4. \qed \enddemo

\demo{Remark}For $H_p(\bold T^n)$ with $0<p<\infty$ we have the
following
situation. When $1<p<\infty$ the orthogonal projection from $L_p(\bold T^n)$
onto $H_p(\bold T^n)$ is bounded so  then $H_p(\bold
T^n)$ is isomorphic to $L_p(\bold T^n)$. This implies in particular that
these spaces are isomorphic to nonatomic lattices.
 When
$0<p<1$ then
$H_p(\bold T^n)$ admit only
purely atomic orders as a $p$-Banach lattices. To see this observe that if $X$
is not a purely atomic $p$-Banach lattice then its Banach envelope
(for definition and properties see \cite{11})  is a Banach lattice which is
not  purely atomic. On the other hand it is  known that
 the Banach envelope of $H_p(\bold T^n)$ is isomorphic to $\ell_1$. For
$n=1$ this can be found in \cite {11} Theorem 3.9, for $n>1$ the proof uses
Theorem 2$'$ of \cite{19} but otherwise is the same; alternatively see
\cite{11} Theorem 3.5, for a proof using bases. When we compare it
with the observation from \cite{1} mentioned in the Introduction, that
$\ell_1$ is not isomorphic to any nonatomic Banach lattice, we
conclude that the spaces $H_p(\bold T^n)$ cannot be isomorphic to any
nonatomic $p$-Banach lattice. \enddemo

\demo{Remark}For the dual spaces $H_1(\bold T^n)^*=BMO(\bold T^n)$ the
situation is rather different.  We first observe the following
proposition: \enddemo
\proclaim{Proposition 3.5}For any Banach space $X$ the spaces
$\ell_1(X)^*(=\ell_{\infty}(X^*))$ and $L_1([0,1],X)^*$ are
isomorphic.\endproclaim

\demo{Proof}Clearly $\ell_1(X)^*$ is isomorphic to a 1-complemented
subspace of $L_1(X)^*.$  Now let $\chi_{n,k}=
\chi_{((k-1)2^{-n},k2^{-n})}$
for $1\le k\le 2^n$ and $n=0,1,\dots.$  Let $T:\ell_1(X)\to L_1(X)$ be
defined by $T((x_n)) =\sum x_n\chi_{m,k}$ where $n=2^m+k-1.$  Let
$L_1(\Cal D_N,X)$ be the subspace of all functions measurable with
respect to the finite algebra generated by the sets
$((k-1)2^{-N},k2^{-N})$ for $1\le k\le 2^N$ and define $S_N:L_1(\Cal
D_N;X)\to
\ell_1(X)$ by
setting $S(x\otimes \chi_{N,k})$ to be the element with $x$ in position
$2^N+k-1$ and zero elsewhere.  Then applying Exercise 7 of II.E of
\cite{22} (cf. \cite{8} Proposition 1), we obtain that $L_1(X)^*$ is
isomorphic to a complemented subspace of $\ell_1(X)^*.$  Then by the
Pelczynski decomposition technique we obtain the proposition.\qed
\enddemo

Now from the Proposition, observe that, since $H_1(\bold T^n)\sim
\ell_1(H_1(\bold T^n))$, we have
$L_1(H_1(\bold T^n))^*
\sim BMO(\bold T^n)$ and clearly this isomorphism induces a nonatomic
(but not order-continuous) lattice structure on $BMO(\bold T^n).$
(It is easy to see that a space which contains a copy of
$\ell_{\infty}$
cannot have an order-continuous lattice structure, because it fails
the separable complementation property.)

\heading 4.   Rad $H_1$ and tent spaces \endheading

The space $H_1(\delta,\ell_2)$ is, as  observed in Section 2, isomorphic
to Rad $H_1$ and has a structure as a nonatomic Banach lattice.
The complex space Rad $H_1$ is easily seen to be isomorphic to the
vector-valued  space  $H_1(\bold T,\ell_2)$ consisting of the boundary
values of the space of all functions $F$  analytic  in the unit disk
$\bold D$ with values in a Hilbert space $\ell_2$ and such that:
$$
\sup_{0<r<1}\int_0^{2\pi}
\|F(re^{i\theta}\|\frac{d\theta}{2\pi}=\|F\|<\infty.$$
To see this isomorphism just note that $H_1(\bold T,\ell_2)$ can be
identified with the space of sequences $(f_n)$ in $H_1$ such that
$$ \|(f_n)\|=\int_0^{2\pi}
(\sum_{n=1}^{\infty}|f_n(e^{i\theta})|^2)^{1/2}\frac{d\theta}{2\pi}<\infty.
$$
This is in turn easily seen to be equivalent to the norm of $\sum r_nf_n$
in $L_2([0,1];H_1)$  (see \cite{16} Theorem 1.d.6).

We now show that a nonatomic Banach lattice isomorphic to Rad $H_1$
arises naturally in
 in harmonic analysis. More precisely we will show that tent space $T^1$
which was introduced and studied by R. Coifman, Y. Meyer and E. Stein
in
\cite{4} is isomorphic to Rad $H_1$. Tent spaces are useful
in some
questions of harmonic analysis (cf. \cite{7} or \cite{21}). They can be
defined over $\bold R^n$ but for the sake of simplicity we will consider them
only over $\bold R$.

Let us fix $\alpha >0$. For $x\in \bold R$ we define
$$\Gamma_\alpha (x) = \{ (y,t) \in \bold R \times \bold R^+ \ :\ |x-y| < \alpha
t \}.$$
Given a function $f(y,t)$ defined on $\bold R \times \bold R^+$ we put
$$\| f\|_\alpha = \int_{\bold R} \big( \int_{\Gamma_\alpha(x)}
|f(y,t)|^2t^{-2} \, dy\, dt \big)^{1/2}\, dx.$$
It was shown in \cite{4} Proposition 4 that for different $\alpha$'s the
norms $\|.\|_\alpha$ are equivalent i.e. for $0<\alpha <\beta <\infty$ there is
a $C=C(\alpha,\beta)$ such that for every $f$ we have
$$\|f\|_\alpha \leq \|f\|_\beta \leq C\|f\|_\alpha . \tag 4.1$$
This implies that the space $T^1 =\{ f(y,t) \ : \ \|f\|_\alpha <\infty \} $
does not depend on $\alpha$. Observe that $T^1$ is clearly a non-atomic Banach
lattice.

Our main result of this Section is
\proclaim{Theorem 4.1 } The space $T^1$ is lattice-isomorphic to
$H_1(\delta, L_2)$ and hence isomorphic to Rad $H_1.$ \endproclaim
Actually for the proof of this Theorem it is natural to work with the dyadic
$H_1$ space on $\bold R$. This space, which we denote
$H_1(\delta_{\infty})$ can be defined as follows:

Let $I_{nk} =[k\cdot 2^n, (k+1)\cdot 2^n]$ for $n,k=0,\pm 1, \pm 2 \dots$ and
let $h_{nk}$ be the function
which is equal to $1$ on the left hand half of $I_{nk}$, $-1$ on the right hand
half of $I_{nk}$ and zero outside $I_{nk}$. In other words $h_{nk}$ is the Haar
system on $\bold R$. The system $\{h_{nk}\}_{n,k =0,\pm 1, \pm 2 \dots}$ is a
complete orthogonal system. For a function $f=\sum_{n,k}a_{nk}h_{nk}$ we define
its $H_1(\delta_{\infty})$-norm by
$$\|f\|= \int_{\bold R} (\sum_{n,k} |a_{nk}|^2 |h_{nk}|^2 )^{1/2}
dt \tag 4.2$$
That this space is isomorphic to the space $H_1(\delta)$ follows from the work
of Sj\"olin and Stromberg \cite{20}.  However, slighlty more is true:

\proclaim{Lemma 4.2}The atomic Banach lattices $H_1(\delta)$ and
$H_1(\delta_{\infty})$ are lattice-isomorphic (or, equivalently
the natural normalized unconditional bases of these spaces are permutatively
equivalent).\endproclaim

\demo{Proof}For any subset $\Cal A$ of $\bold Z^2$ write $H_{\Cal A}$ for
the closed linear span of $\{h_{nk}:(n,k)\in\Cal A\}$ in
$H_1(\delta_{\infty}).$  For $m\in\bold Z$ let $\Cal
A_m=\{(n,k):I_{nk}\subset [2^{-m-1},2^{-m}]\}$ and $\Cal B_m=\{(n,k):I_{nk}
\subset [-2^{-m},-2^{-m-1}].$  Let $\Cal D=\cup_{m\in\bold Z}(\Cal
A_m\cup\Cal
B_m)$ and $\Cal D_+=\cup_{m\ge 0}\Cal A_m.$  Then it is clear that $H_{\Cal
D}$ and $H_{\Cal D_+}$ are each lattice isomorphic to $\ell_1(H_1(\delta)).$
Now $H_1(\delta_{\infty})$ is lattice isomorphic to $H_{\Cal D}\oplus
H_{\Cal E}$ where $\Cal E=\{(m,0),\ (m,-1):m\in\bold Z\}$.  It is easy
to show that $H_{\Cal E}$ is lattice isomorphic to $\ell_1.$  Similarly
$H_1(\delta)$ is lattice-isomorphic to $H_1(\Cal D_+)\oplus\ell_1$ and this
completes the proof of the lemma.\qed \enddemo

\demo{Remark}Note also that $H_1(\delta)$ is lattice-isomorphic to
$\ell_1(H_1(\delta)).$  \enddemo

   \demo{Proof of the Theorem} We will prove that $T^1$ is
lattice-isomorphic to $H_1(\delta_{\infty},L_2)$.  Let us introduce
squares
$A_{nk}
\subset
\bold R
\times \bold R^+$ defined as $A_{nk} = I_{nk} \times [2^n,2^{n+1}]$ for $n,k
=0,\pm 1, \pm 2, \dots$. It is geometrically clear that squares
$\{A_{nk}\}_{n,k =0, \pm 1 ,\pm 2, \dots}$ are essentially disjoint and that
they cover $\bold R\times \bold R^+$. For $j=0,1,2$ we define
$$A^j_{nk} = [(k+ \tfrac j3)2^n, (k+\tfrac{j+1}3 )2^n] \times [2^n, 2^{n+1}].$$
Note that in this way we divide each $A_{nk}$ into three essentially
disjoint
rectangles. Let $D^j = \bigcup_{n,k} A^j_{nk}$. Let $T^1_j$ be the subspace of
$T^1$ consisting of all functions whose support is contained in $D^j$.  Clearly
$T^1 =T^1_0 \oplus T^1_1 \oplus T^1_2$, so it is enough to show that
$T^1_j$ is lattice-isomorphic to $H_1(\delta_{\infty},L_2).$

We write $f^j \in T^1_j$  as $f^j =\sum_{n,k} f^j_{nk}$ where $f^j_{nk}=f^j
\cdot \chi _{A^j_{nk}}$. We start with $j=1$. For any $\alpha >0$ we have
$$\split
\|f^1\|_{\alpha} & =\int_{\bold R}(\int_{\Gamma_\alpha(x)} |f^1(y,t)|^2
t^{-2}dy\, dt)^{1/2} dx\\
&=\int_{\bold R}(\int_{\Gamma_\alpha(x)} \sum_{n,k}|f^1_{nk}(y,t)|^2 t^{-2}
dy\, dt)^{1/2}dx\\
& = \int_{\bold R} ( \sum_{nk} \int_{\Gamma_\alpha(x)} | f^1_{nk}(y,t)|^2
t^{-2}dy\,dt)^{1/2} dx. \endsplit \tag 4.3$$
If we now take $\alpha =\frac23$ we have $\Gamma_\alpha(x) \supset A^1_{nk}$
for all $x \in I_{nk}$, so from \thetag{4.3} we get
$$\|f^1\|_\alpha \geq \int_{\bold R} ( \sum_{nk} \chi_{I_{nk}}(x)
\int_{A^1_{nk}}|f^1_{nk}(y,t)|^2 t^{-2} dy\, dt )^{1/2} dx. \tag 4.4$$
On the other hand when we take $\alpha =\frac16$ we have $\Gamma_\alpha(x) \cap
A^1_{nk} = \emptyset $ for all $x \notin I_{nk}$, so from \thetag{4.3} we get
$$\|f^1\|_\alpha \leq \int_{\bold R} ( \sum_{n,k} \chi_{I_{nk}}(x)
\int_{A^1_{nk}}
| f^1_{nk} (y,t) |^2 t^{-2} dy\, dt)^{1/2} dx. \tag 4.5$$
For each $(n,k)$ the subspace of $T^1$ consisting of functions supported on
$A^1_{nk}$ is easily seen to be isometric to the Hilbert space. If we fix an
isometry between this space and $\ell_2$ we obtain from \thetag{4.2},
\thetag{4.3} and \thetag{4.4} that $T^1_1$ is lattice-isomorphic to
$H_1(\delta_{\infty},L_2)$. In
order to complete the proof of the Theorem it is enough to show that $T^1_0$
and $T^1_2$ are lattice-isomorphic to $T^1_1$. This isomorphism can be
given by
$\sum_{nk} f^j_{nk} \mapsto \sum_{nk} f^1_{nk}$. The fact that this map is
really an isomorphism follows from: \enddemo
\proclaim{Lemma 4.3} Let $\phi (t)$ be a uniformly bounded measurable
function
on $\bold R^+$. For a function $f$ defined on $\bold R \times \bold R^+$ we
define
$$A_\phi(f)(y,t) = f(y+ t\phi(t),t).$$
Then $A_\phi \,: T^1 \longrightarrow T^1$ is a continuous linear operator.
\endproclaim
\demo{Proof of the Lemma} Since
$$\split
&\int_{\Gamma_\alpha(x)} | A_\phi (f)(y,t) |^2 t^{-2} dy\, dt \\
 &= \int_{\bold R^+} (t^{-2} \int_{x-\alpha t}^{x+\alpha t} | A_\phi
(f)(y,t)|^2 dy)dt\\
& = \int_{\bold R^+} (t^{-2} \int_{x- \alpha t - t \phi(t)}^{x+\alpha t - t
\phi (t)}| f(y,t)|^2 dy)dt \\
&\leq \int_{\bold R^+} (t^{-2} \int _{x-(\|\phi \|_\infty +\alpha
)t}^{x+(\|\phi \|_\infty +\alpha)t} | f(y,t)|^2 dy)dt \\
&= \int_{\Gamma_{\alpha+\|\phi\|_\infty}(x)} |f(y,t)|^2 t^{-2} dy\, dt
\endsplit $$
the Lemma follows.\qed \enddemo

\enddocument